\documentclass[reqno, 11pt]{article}
\usepackage[centertags,tbtags]{amsmath}
\usepackage{amssymb} 
\usepackage{amsthm}
\usepackage{txfonts}
\usepackage{fullpage}
\usepackage{enumitem}
\usepackage{bbm} 
\usepackage[usenames,dvipsnames]{color}
\usepackage[pdftex]{thumbpdf}       
\usepackage{ifpdf}

\ifpdf
\usepackage[pdftex]{hyperref}
\hypersetup{%
pdftitle={Regenerative Processes},
pdfauthor={Maria Vlasiou},
pdfkeywords={EORMS},
bookmarks=true,%
bookmarksnumbered=true,
pdfstartview={FitH},
colorlinks=true,
citecolor=Plum,
linkcolor=RedOrange,
urlcolor=RawSienna,
bookmarksopen=true,%
bookmarksopenlevel=1,     
unicode=true,%
breaklinks=true,%
}%
\else
\newcommand\phantomsection\relax
\newcommand{\url}[1]{#1}
\newcommand{\href}[2]{#2}
\usepackage{nohyperref}
\fi

\theoremstyle{plain}              
\newtheorem{theorem}{Theorem}

\newtheorem{proposition}{Proposition}

\theoremstyle{definition}
\newtheorem{remark}{Remark}
\newtheorem{ex}{Example}
\newtheorem{defi}{Definition}

\newcommand{\e}{\mathbb{E}}
\newcommand{\p}{\mathbb{P}}
\renewcommand{\d}{\,\mathrm{d}}

\begin{document}
\title{Regenerative Processes}
\author{Maria Vlasiou\thanks{Dept.\ of Mathematics \& Computer Science, Eindhoven University of Technology, P.O.\ Box 513, 5600 MB Eindhoven, The Netherlands, \href{mailto:m.vlasiou@tue.nl}{m.vlasiou@tue.nl}}}
\date{\today}
\maketitle

\begin{abstract}
We review the theory of regenerative processes, which are processes that can be intuitively seen as comprising of i.i.d.\ cycles. Although we focus on the classical definition, we present a more general definition that allows for some form of dependence between two adjacent cycles, and mention two further extensions of the second definition. We mention the connection of regenerative processes to the single-server queue, to multi-server queues and more generally to Harris ergodic Markov chains and processes. In the main theorem, we pay some attention to the conditions under which a limiting distribution exists and provide references that should serve as a starting point for the interested reader.
\end{abstract}

\phantomsection
\addcontentsline{toc}{section}{Classical regenerative processes}
\section*{Classical regenerative processes}
A stochastic process $\{X(t),\, t\geqslant 0\}$ is intuitively a regenerative process if it can be split into i.i.d.\ cycles. That is, we assume that a collection of time points exists, so that between any two consecutive time points in this sequence, (i.e.\ during a cycle), the process $\{X(t),\, t\geqslant 0\}$  has the same probabilistic behaviour. There are several proposed definitions which attempt to make this description precise; see Asmussen~\cite{asmussen-APQ}, \c{C}inlar~\cite{cinlar-ISP} or Thorisson~\cite{thorisson83}. The classical definition of regenerative processes, which is given below, also demands that the process is independent of its history up to the previous regeneration point, and of the choice of the point itself. This assumption has been relaxed significantly and the notion of regenerative processes that we understand nowadays is wider; we will sketch these ideas in the following sections. This presentation is influenced by many books, which can be found in the references, and draws heavily from Sigman and Wolff~\cite{sigman93}.

The power of the concept of regenerative processes lies in the existence of a limiting distribution under conditions that are mild and usually easy to verify. For example, in continuous time it is only required that the mean cycle length of the embedded renewal process is finite, that the cycle length distribution is non-lattice and that the sample paths of the regenerative process satisfy some conditions (e.g.\ right-continuity with left-hand limits) that are automatic in most examples; see Asmussen \cite[Chapter VI]{asmussen-APQ}.

Classic examples of regenerative processes are ergodic Markov chains and processes. For example, take a Markov chain or a Markov process with a countable state space $E$ and fix a state $i$. Assume that the process started at state $i$. Then every time at which state $i$ is entered is a time of regeneration for the Markov process. As a second example, consider a single-server queuing system with Poisson arrivals and i.i.d.\ service times. Suppose that we start observing the system at a departure instant which left behind $i$ customers. Then every time a customer departs leaving behind him $i$ customers is a regeneration epoch and the process describing the queue length after such a time has exactly the same probability law as the one it had at time zero.

We start with the classical definition of a regenerative process.

\begin{defi}[Classical definition]\label{def:1}
A regenerative process $\{X(t), \,t\geqslant 0\}$ with state space $E$ (endowed with a $\sigma$-field $\mathcal{E}$) is a stochastic process on an underlying probability space $(\Omega, \mathcal{F}, \p)$ with the following properties: there exists a random variable $R_1>0$ such that
\begin{enumerate}[topsep=0pt, label=(\roman*), itemsep=0pt]
\item $\{X(t+R_1), \,t\geqslant 0\}$  is independent of $\{X(t), \,t\leqslant R_1\}$  and $R_1$, \label{pr:1}
\item $\{X(t+R_1),\, t\geqslant 0\}$  is stochastically equivalent to $\{X(t),\, t\geqslant 0\}$.\label{pr:2}
\end{enumerate}
\end{defi}
\textsl{Stochastically equivalent} means that the two stochastic processes in \ref{pr:2} have the same joint distributions. Because of property \ref{pr:1}, $R_1$ is a stopping time for $\{X(t), \,t\geqslant 0\}$. For a rigorous definition of property \ref{pr:1} see also Sigman and Wolff~\cite{sigman93}. A standard assumption that is often made is that the paths of $X(t)$ are right-continuous with left-hand limits almost surely. This ensures that the sample paths are continuous except possibly on a set of Lebesgue measure 0. We will comment on the necessity of this assumption in the main theorem later on.

The random variable $R_1$ is called a \textsl{regeneration epoch}, and its existence implies the existence of a sequence of regeneration epochs $\{R_n\}$. This result is carefully proven in Sigman, Thorisson and Wolff~\cite{sigman94} and Svertchkov~\cite{svertchkov93}. We say that the process \textsl{regenerates} at such epochs. In other words, $\{X(t), \, t\geqslant 0\}$ and $\{X(t+R_n),\, t\geqslant 0\}$ are stochastically identical, and due to property~\ref{pr:1}, $\{X(t+R_n),\, t\geqslant 0\}$ is independent of $\{X(t),\, 0 \leqslant t\leqslant R_n\}$. The regeneration epochs generate a renewal process $S_n=R_1+\cdots+R_n$ (see Section 2.1.4.1 of EORMS), which is called the \textsl{embedded renewal process} for the regenerative process $\{X(t), \, t\geqslant 0\}$. Moreover, the interval $[S_{n-1},S_n)$ is called the \textsl{$n$-th regenerative cycle}, where we have assumed that $S_0=0$. A regenerative process with finite mean cycle length is called \textsl{positive recurrent}.

It will sometimes be convenient to allow the process during the first cycle to be different; this simply means that the process during the first cycle may have a distribution that is different from that of subsequent cycles. We call such a process \textsl{delayed regenerative}, if the emphasis on the assumptions of the first cycle is needed. From a time-average point of view, the first cycle is not important, provided that it eventually ends. That is, we assume that $R_1$ is a proper random variable. There are no new tools needed; in handling a delayed regenerative process first we condition the event or expectation in question on the time $R_1$ of the first regeneration, and then we use the fact that at $R_1$ there starts an ordinary regenerative process. Thus, we use the term \textit{regenerative} to cover both cases of delayed or not regenerative processes.

\begin{remark}
If the distribution of the process during the first special cycle and the length of that cycle are constructed in a particular way, then we can construct a stationary regenerative process, which implies that the distribution of $\{X(t), \, t\geqslant 0\}$ is independent of $t$ for all $t$; see Wolff~\cite{wolff-SMTQ} for a construction.
\end{remark}

An important fact that follows directly from the classical definition is that functions of regenerative processes are also themselves regenerative.

\begin{proposition}[Inheritance of regenerations]
If $\{X(t), \,t\geqslant 0\}$ with state space $E$  is a regenerative process with regeneration epochs $R_n$, then $\{f(X(t)), \,t\geqslant 0\}$ is also a regenerative process (with some state space $E^\prime$) with the same regeneration epochs for any $f:E\to E^\prime$.
\end{proposition}

For instance, take $\{X(t), \,t\geqslant 0\}$ to be a regenerative process and consider the function $\hat{X}(t)=\mathbbm{1}(X(t)\in B)$ for some set $B$. Since $X(t)$ is regenerative, then so is $\hat{X}(t)$, which can help one compute the distribution of one regenerative process as the expectation of another, i.e.\ $\p(X(t)\in B)=\e[\hat{X}(t)]$. This proposition illustrates a more general observation that the embedded renewal times of a regenerative process need not be stopping times with respect to the evolution of $X(t)$. Observe that in contrast a function of a Markov process is usually not a Markov process itself. Nevertheless, it remains regenerative if the Markov process is so.

Other than computing the limiting distribution of one regenerative process via another, one could use the inheritance of regenerations to add a cost or reward function to the process. Namely, if we assume that the system earns rewards (that may be negative) at a rate $c(x)$ whenever it is in state $x$, then the total reward up to time $t$ is given by $\int_0^t c(X(y)) \d y$, where $\{c(X(y)), y\geqslant0\}$ is a regenerative process. For finite mean cycle lengths, the time average of $c(x)$ defined to be equal to $\lim_{t\to\infty}\int_0^t c(X(y)) \d y/t$ can be found by the renewal-reward theorem (see Theorem 1 in Section 2.1.4.6 of EORMS). In particular, it can be shown with standard methods that
$$
\lim_{t\to\infty}\int_0^t c(X(y)) \d y/t = \int_{-\infty}^\infty c(y) \d F(y),
$$
where $F(x)$ is the limiting distribution of the regenerative process $\{X(t), \,t\geqslant 0\}$, which is given below in Theorem~\ref{th:1} by choosing $(-\infty, x]$ for the set $B$ appearing in Equation~\eqref{eq:1}.

The primary use of the key renewal theorem (see Section 2.1.4.3 of EORMS) is in characterising the limiting behaviour of regenerative processes. Renewal theory is the main tool for studying regenerative processes in the absence of further properties. Conversely, it is this applicability to regenerative processes which makes renewal theory one of the most important tools in elementary probability theory. The time average properties of regenerative processes are straightforward consequences of the renewal reward theorem (which is based on the strong law of large numbers). This follows from the fact that the cycles, and the process during each cycle, are i.i.d. We proceed with the main theorem in this section, which gives the limiting distribution of a regenerative process.

\begin{theorem}[Limiting behaviour of regenerative processes]\label{th:1}
Suppose that the process $\{X(t), \,t\geqslant 0\}$ is regenerative with regeneration epochs with finite mean, i.e.\ $\e[R_1]<\infty$. We then have
\begin{enumerate}[label=(\alph*)]
\item
The limiting distribution of the regenerative process in the time average sense exists; i.e.\ the limit $\lim_{t\to\infty} \frac{1}{t} \int_0^t \mathbbm{1} (X(s)\in B) \d s$ exists for all measurable sets $B$ , $B\in\mathcal{E}$, where $ \mathbbm{1}(\cdot)$ is the indicator function. For $B=(-\infty,x]$ we denote it by $F(x)$. Moreover, if we denote by $X_\infty$ a generic random variable with this limiting distribution, we have that
 \begin{equation}\label{eq:1}
    \p(X_\infty\in B)= \frac{1}{\e[R_1]}\e[\int_0^{R_1} \mathbbm{1}(X(s)\in B) \d s, \quad B\in\mathcal{E}.
\end{equation}

\item \label{item:b}
If the state space $E$ is a complete separable metric space and the process $\{X(t), \,t\geqslant 0\}$ has right-continuous paths, then in case the cycle length distribution is non-lattice, we have that $\p(X(t)\in B)$ converges weakly to $\p(X_\infty\in B)$. Moreover, for a measurable function $f:E\to\mathbbm{R}$ we have that
 \begin{equation}\label{eq:2}
\lim_{t\to\infty} \e[f(X(t))]=\frac{1}{\e[R_1]}\e[\int_0^{R_1} f(X(s)) \d s,
\end{equation}
which for the special case where $f$ is the indicator function, yields that $\lim_{t\to\infty} \p(X(t)\in B)=\p(X_\infty\in B)$ is equal to the right hand side of \eqref{eq:1}.

\item
If the cycle length distribution is spread-out, then $\p(X(t)\in B)$ converges  to $\p(X_\infty\in B)$ in total variation (i.e.\ the convergence is uniform over all Borel sets $B$).
\end{enumerate}
\end{theorem}
For a regenerative process in discrete time, the above theorem holds in case the cycle length is aperiodic. If it is periodic with period $d$, the above limit holds if $t=nd$ and $n\to\infty$, see also Asmussen~\cite{asmussen-APQ} and Kulkarni~\cite{kulkarni-MASS}. Serfozo~\cite{serfozo-BASP} derives \eqref{eq:2} without the assumptions under \ref{item:b}, but with the assumptions that the cycle length distribution is non-lattice and that for $M=\sup\{\,|f(X(t))|\,:t\leqslant R_1\}$, $M$ and $MR_1$ have finite means.

From this theorem we see that the limiting distribution for a regenerative process, $\p(X_\infty \in B)$, is equal to the expected time that the process spent in $B$ over the expected inter-renewal time. Evidently, if the embedded renewal process is recurrent but with $\e[R_1]=\infty$, then we have that $\p(X_\infty \in B)=0$. One could use the above theorem to derive the Laplace-Stieltjes transform of the limiting distribution of $\{X(t), \,t\geqslant 0\}$ by the appropriate choice of the function $f$ appearing in \ref{item:b}.

This theorem is also true, with slight modifications, for delayed regenerative processes. We namely have to integrate over a proper regenerative cycle, say $[S_1,S_2)$, rather than the first special cycle.

\begin{ex}[Age process]
Evidently, the process $Z_t$ giving the time from $t$ until the next renewal in a renewal process is a regenerative process.
\end{ex}

\begin{ex}[GI/GI/1 queue]
Let $\{X(t), \,t\geqslant 0\}$ be the number of customers in a GI/GI/1 queue and let $S_n$ be the $n$-th time when a customers enters an empty system. Provided that the occupation rate of the system is less than 1 (i.e.\ that the mean service time is less than the mean interarrival time), the mean cycle length is finite. If the process starts with a customer entering an empty system at time 0, then $\{X(t), \,t\geqslant 0\}$ is a regenerative process; otherwise it is a delayed regenerative process. The state space of this regenerative process is given by $\{0,1,2,\ldots\}$. Then, under the assumptions of the theorem above
$$
\lim_{t\to\infty} \p(X(t)=j)= \frac{\e[\mbox{time in state $j$ during a cycle}]}{\e[\mbox{time of a cycle}]}=\lim_{t\to\infty} \frac{\e[\mbox{time in $j$ during $(0,t)$}]}{t}.
$$
That is, $\lim_{t\to\infty} \p(X(t)=j)$ is not only the probability that the regenerative process is in state $j$ in steady-state (i.e.\ the probability that there are $j$ customers in the system in steady-state), but also the long-run proportion of time that the process was in state $j$. To prove this result, one simply needs to assume that a reward is earned at rate 1 each time the process is in state $j$, which generated a renewal reward process and the equality is then a consequence of the renewal-reward theorem (or of Theorem~\ref{th:1} above, for a particular choice of $f$).
\end{ex}

\begin{ex}[I.i.d.\ random variables]
Any independent and identically distributed sequence of random variables is a regenerative process. The embedded renewal process is simply the counting process of the sequence ($S_n=n$). If these random variables are continuous, we observe that the process returns to a fixed state with probability zero, and thus, the regeneration epochs cannot be defined as returns to any particular fixed state; contrast this with the fact that in the GI/GI/1 queue example above, returns to the fixed state ``the number of customers in the queue is zero'' were chosen to be the regenerative epochs.
\end{ex}

As mentioned earlier, one can find several examples of regenerative processes in Markov chains and Markov processes.

\begin{remark}
There is a version of the central limit theorem for regenerative processes that can be used to construct confidence intervals when simulating regenerative processes; see Sigman and Wolff~\cite{sigman93} and Wolff~\cite{wolff-SMTQ} for details.
\end{remark}

\phantomsection
\addcontentsline{toc}{section}{One-dependent regenerative processes}
\section*{One-dependent regenerative processes}
The definition of regenerative processes has been extended to allow for one-dependent cycles, which connects regenerative processes to Harris-recurrent Markov chains. Many-server queues formed the main motivation to extend the classical notion of regeneration to this definition \cite{sigman93}. Previously, for the GI/GI/1 queue we had taken the times when the system empties as the regeneration epochs. Suppose now you have a GI/GI/2 queue which is stable, i.e., the offered load to the system per server is less than 1. Take a completely deterministic system where service takes 1.5 time units and customers arrive every 1 time unit. Let the first arrival occur at time 0. Then we observe that this system never empties. Thus, we are not able to define the same regeneration epochs (i.e.\ an empty system) as for the single-server queue. The question is now if one should infer that the queue length process is not regenerative.

The definition of regenerative processes we have used so far demands that $\{X(t+R_1), \,t\geqslant 0\}$ is independent both of $R_1$ and of $\{X(t), \,t\leqslant R_1\}$ (see Definition~\ref{def:1}). A more general definition is given in Asmussen~\cite{asmussen-APQ} where the requirement now is that $\{X(t+R_1), \,t\geqslant 0\}$ is independent only of $R_1$ but not necessarily of $\{X(t), \,t\leqslant R_1\}$. Under both definitions, $\{X(t+R_1), \,t\geqslant 0\}$ is stochastically equivalent to $\{X(t), \,t\geqslant 0\}$. As a result, the embedded renewal process has i.i.d.\ cycles, but the regenerative process might be dependent. Asmussen~\cite{asmussen-APQ} shows that under this definition, the process is \textsl{one-dependent}, namely, the process during adjacent cycles is dependent, but during non-adjacent cycles it is independent, see also Nummelin~\cite{nummelin94}. If a process is regenerative under Definition~\ref{def:1}, then it is naturally also so under this definition, while evidently, the converse is not true.

\begin{defi}[One-dependent definition]\label{def:2}
A regenerative process $\{X(t), \,t\geqslant 0\}$ with state space $E$ (endowed with a $\sigma$-field $\mathcal{E}$) is a stochastic process on an underlying probability space $(\Omega, \mathcal{F}, \p)$ with the following properties: there exists a random variable $R_1>0$ such that
\begin{enumerate}[topsep=0pt, label=(\roman*), itemsep=0pt]
\item $\{X(t+R_1), \,t\geqslant 0\}$  is independent of $R_1$\label{pr:2a}
\item $\{X(t+R_1),\, t\geqslant 0\}$  is stochastically equivalent to $\{X(t),\, t\geqslant 0\}$.
\end{enumerate}
\end{defi}

Since the cycle lengths $R_n$ are still i.i.d., they generate an embedded renewal process. Thus, \eqref{eq:2} still holds; the regenerative process converges in distribution to its limit (and also in total variation). This definition can be used to characterise Harris recurrent Markov chains, as a Markov chain is positive Harris recurrent if and only if it is regenerative under the one-dependent definition and ergodic for every initial state $x$, with the distribution of the process during a cycle (possibly excluding the first one for the delayed case) not depending on the initial state $x$; see~\cite{sigman90,sigman93}. We usually think of queues as continuous-time processes, but we can embed discrete-time processes in them, as for example instances when a customer arrives. Embedded processes for several queuing models turn out to be regenerative under the one-dependent definition, and with finite mean cycle lengths \cite{sigman88a,sigman88}. Returning to the GI/GI/2 example mentioned above, it has been shown \cite{charlot78} that the vector giving the workload for each server is a Harris ergodic Markov chain and thus it has a regenerative structure; for details see \cite{sigman88,wolff-SMTQ}.

\phantomsection
\addcontentsline{toc}{section}{Further reading}
\section*{Further reading}
The interested reader should start with the comprehensive and in depth review of regenerative processes by Sigman and Wolff~\cite{sigman93}; several important references on the topic can be found there -- we present only a few of them. Smith~\cite{smith55,smith58} constructs a regenerative process under the classical definition using random tours; see also \cite{resnick-ASP}.  As this topic is an integral part of stochastic processes, most books on stochastic processes have some treatment of (classical) regenerative processes; see for example \cite{cinlar-ISP,cohen-RPQT,kulkarni-MASS,resnick-ASP,ross-SP,serfozo-BASP}. References on how to simulate a stationary version of a regenerative process have been omitted; few of them can be found to the sources cited here.

The notion of regeneration has been extended to even more general definitions. Rather than demanding that the regenerative epochs $R_n$ form an i.i.d.\ sequence (cf.\ property \ref{pr:2a} in Definition~\ref{def:2}), one can extend the notion of the dependence of $\{X(t)\}$ on the process in adjacent cycles, to one-dependence of the cycle lengths themselves. In Sigman~\cite{sigman90} it is shown that these processes are related to \textsl{continuous} time Harris-recurrent Markov processes; (compare this to what has been mentioned in the previous section that Definition~\ref{def:2} is related to Harris ergodic chains, i.e.\ in discrete time). Further examples of this extended notion of one-dependence of regenerative processes include the superposition of two independent renewal processes where at least one of them has a spread-out cycle length distribution.

Another way of describing the regenerative structure for a recursive sequence as this was given under Definition~\ref{def:2} is by the \textsl{renewing events} or \textsl{renovation method} \cite{borovkov-AMQT,foss91}. This amounts to first \emph{algebraically} describing the independence of the past property rather than doing so in a probabilistic fashion; see \cite{sigman93} for an example.

A further extension is the notion of \textsl{synchronous processes}; a stochastic process $\{X(t), \,t\geqslant 0\}$ is called synchronous if there is a random variable $R_1>0$ such that $\{X(t+R_1), \,t\geqslant 0\}$ is stochastically equivalent to $\{X(t), \,t\geqslant 0\}$, possibly by enlarging the probability space if necessary. That is, we no longer require property \ref{pr:2a} in Definition~\ref{def:2}; see \cite{franken-QPP,glynn92,rolski-SRPAPP}.

Other notions of regenerative processes worth mentioning are those of Kingman~\cite{kingman-RP,kingman06} in which the kind of regenerative structure found in discrete state Markov processes is generalised by allowing any of the continuum of times during the visit to a recurrent state to serve as a regeneration epoch and that of Thorisson~\cite{thorisson83} who introduces the notion of a time-inhomogeneous regenerative process, which can apply for example to an inhomogeneous Markov chain with a recurrent state.

The list of references here is by no means exhaustive, but should serve only as a starting point.

\phantomsection
\addcontentsline{toc}{section}{References}

\end{document}